\documentclass{article}
\usepackage{stmaryrd}
\usepackage{mathrsfs}
\usepackage[centertags]{amsmath}
\usepackage{amsfonts, dsfont}
\usepackage{amssymb}
\usepackage{amsthm}

\input xy
\xyoption{all}

\theoremstyle{plain}
\newtheorem{thm}{Theorem}[section]
\newtheorem{cor}[thm]{Corollary}
\newtheorem{lem}[thm]{Lemma}
\newtheorem{prop}[thm]{Proposition}

\theoremstyle{definition}
\newtheorem{defn}[thm]{Definition}
\newtheorem{remark}[thm]{Remark}
\newtheorem*{ack}{Acknowledgments}

\newcommand{\bd}{\begin{defn}}
\newcommand{\ed}{\end{defn}}
\newcommand{\bl}{\begin{lem}}
\newcommand{\el}{\end{lem}}
\newcommand{\bp}{\begin{prop}}
\newcommand{\ep}{\end{prop}}
\newcommand{\bt}{\begin{thm}}
\newcommand{\et}{\end{thm}}
\newcommand{\bc}{\begin{cor}}
\newcommand{\ec}{\end{cor}}
\newcommand{\br}{\begin{remark}}
\newcommand{\er}{\end{remark}}
\newcommand{\bdi}{\begin{diagram}}
\newcommand{\edi}{\end{diagram}}
\newcommand{\beq}{\begin{equation}}
\newcommand{\eeq}{\end{equation}}
\newcommand{\ba}{\begin{array}}
\newcommand{\ea}{\end{array}}
\newcommand{\bpf}{\begin{proof}}
\newcommand{\epf}{\end{proof}}

\newcommand{\Z}{\mathds{Z}}
\newcommand{\Q}{\mathds{Q}}
\newcommand{\Zp}{\mathds{Z}_{p}}
\newcommand{\Qp}{\mathds{Q}_{p}}
\newcommand{\al}{\alpha}
\newcommand{\Ga}{\Gamma}
\newcommand{\ga}{\gamma}
\newcommand{\La}{\Lambda}
\newcommand{\la}{\lambda}

\newcommand{\M}{\mathfrak{M}}

 \DeclareMathOperator{\Gal}{Gal}
\DeclareMathOperator{\Hom}{Hom} \DeclareMathOperator{\rank}{rank}
\DeclareMathOperator{\corank}{corank}
\DeclareMathOperator{\Ext}{Ext} \DeclareMathOperator{\Ann}{Ann}

\newcommand{\ot}{\otimes}
\newcommand{\ilim}{\displaystyle \mathop{\varinjlim}\limits}

\newcommand{\cyc}{\mathrm{cyc}}
\newcommand{\cts}{\mathrm{cts}}

\newcommand{\lra}{\longrightarrow}

\newcommand{\ps}[1]{\llbracket #1 \rrbracket}

\begin{document}

\title{On completely faithful Selmer groups of elliptic curves and Hida deformations}
\author{Meng Fai Lim\footnote{Department of Mathematics, University of Toronto, 40 St. George St.,
Toronto, Ontario, Canada M5S 2E4. E-mail: mengfai.lim@utoronto.ca}}
\date{}
\maketitle

\begin{abstract} \footnotesize
\noindent
 In this paper, we study completely faithful torsion $\Zp\ps{G}$-modules
with applications to the study of Selmer groups. Namely, if $G$ is a
nonabelian group belonging to certain classes of polycyclic pro-$p$
group, we establish the abundance of faithful torsion
$\Zp\ps{G}$-modules, i.e., non-trivial torsion modules whose global
annihilator ideal is zero. We then show that such
$\Zp\ps{G}$-modules occur naturally in arithmetic, namely in the
form of Selmer groups of elliptic curves and Selmer groups of Hida
deformations. It is interesting to note that faithful Selmer groups
of Hida deformations do not seem to appear in literature before. We
will also show that faithful Selmer groups have various arithmetic
properties. Namely, we show that faithfulness is an isogeny
invariant, and we will prove ``control theorem" results on the
faithfulness of Selmer groups over a general strongly admissible
$p$-adic Lie extension.



\end{abstract}

\section{Introduction}

Throughout the paper, $p$ will always denote an odd prime. Let $E$
be an elliptic curve defined over $\Q$ which has good ordinary
reduction at the prime $p$. The Iwasawa main conjecture predicts
that the Mazur-Swinnerton-Dyer $p$-adic $L$-function
$\mathcal{L}_p(E)$ associated to $E$ can be interpreted as an
element of the Iwasawa-algebra $\Zp\ps{\Gal(\Q^{\cyc}/\Q)}$ of the
cyclotomic $\Zp$-extension $\Q^{\cyc}$ of $\Q$ and is a generator of
the characteristic ideal of the Pontryagin dual $X(E/\Q^{\cyc})$ of
the Selmer group of $E$ over $\Q^{\cyc}$ (see \cite{MazSw}).
Furthermore, if $X(E/\Q^{\cyc})$ does not have any nonzero
pseudo-null submodule, it will follow from the main conjecture that
the $p$-adic $L$-function $\mathcal{L}_p(E)$ annihilates
$X(E/\Q^{\cyc})$. We should mention that the Iwasawa main conjecture
in this context has been well understood and largely proven (see
\cite{Ka04, R, SU}).

It is natural to consider generalization of the above by considering
field extensions $F_{\infty}$ of some number field $F$ whose Galois
group $G = \Gal(F_{\infty}/F)$ is a nonabelian $p$-adic Lie group,
and this has been the central theme in noncommutative Iwasawa
theory. One of the earliest approach towards understanding and
formulating this theory is to investigate the global annihilator
ideal of $X(E/F_{\infty})$ (for instance, see \cite{CSS, Ha}).
Inspired by the cyclotomic situation, it was then hoped that such an
investigation might give some insight to the noncommutative $p$-adic
$L$-function which is, even today, still largely conjectural in most
situations (although one now has a slightly better understanding of
the shape of the $p$-adic $L$-functions and the form of the main
conjecture via an algebraic $K$-theoretical approach; see \cite{BV,
CFKSV, FK, Ka05}). As it turns out, such an approach via global
annihilators had been shown to be \textit{not} feasible in general.
In fact, Venjakob was able to establish the existence of a class of
modules over the Iwasawa algebra of the nonabelian group
$\Zp\rtimes\Zp$ which \textit{cannot} be annihilated by a single
global element in the Iwasawa algebra (see \cite{V03}). Building on
this work, he and Hachimori were able to give examples of dual
Selmer groups of elliptic curves over a false Tate extension which
do not have a nonzero global annihilator (see \cite{HV}).

In this paper, following the footsteps of Venjakob, we will
establish the nonexistence of global annihilators for a class of
modules over the Iwasawa algebra of a nonabelian group which is an
extension of a polycyclic pro-$p$ group by $\Zp\rtimes\Zp$ (see
Theorem \ref{faithful modules}). The technique we used in
establishing this result derives essentially from \cite{V03}. In
fact, to prove our main result, we will also require the case of
$\Zp\rtimes\Zp$ which was first established by Venjakob in
\cite{V03}. Then as in the paper of Hachimori and Venjakob
\cite{HV}, we apply our result to obtain examples of completely
faithful Selmer groups of elliptic curves and Hida deformations over
noncommutative $p$-adic Lie extensions (see Theorems \ref{faithful
Selmer groups} and \ref{faithful Selmer groups for Hida}). To the
best of the author's knowledge, completely faithful Selmer groups of
Hida deformations do not seem to be observed in literature before.
We mention that we can also find examples of Selmer groups of
elliptic curves which are faithful but not completely faithful. We
also mention that our results can be applied to obtain completely
faithful Selmer groups for $p$-adic representations defined over
coefficient rings $\Zp\ps{X_1, X_2, ..., X_n}$. However, in this
paper, we will content ourselves mainly with Hida deformations and a
short remark in the general aspect.

For the remainder of the paper, we discuss further properties on the
faithfulness of Selmer groups. In particular, we show that
faithfulness is an isogeny invariant (see Proposition \ref{faithful
isogeny}). On the other hand,  we will give an example to show that
completely faithfulness is not an isogeny invariant. In the final
section of the paper, we will prove some ``control theorem" type
results for the faithfulness of Selmer groups (see Propositions
\ref{faithful control} and \ref{faithful control for Hida}). It
would seem that ``control theorem" type results for faithfulness of
Selmer group have not been observed in literature before.

We should also mention that completely faithful modules and Selmer
groups of elliptic curves over Iwasawa algebras of compact $p$-adic
Lie group other than the ones considered in this paper have also
been studied in \cite{A, BZ}. Our results here may therefore be
viewed as complement to the results there.

We end the introductory section discussing some (negative)
consequences and significance of our results. Let $E$ be an elliptic
curve defined over $\Q$ which has good ordinary reduction at the
prime $p$ and set $F=\Q(\mu_p)$. Let $F_{\infty}$ be a false Tate
extension of $F$ in the sense of \cite{HV}. As shown loc. cit.,
there are cases of $X(E/F_{\infty})$ being completely faithful. One
may then naively consider adjoining multiple $\Zp^r$-extensions of
$F$ to $F_{\infty}$ and perhaps hope to obtain nontrivial global
annihilator of the Selmer groups which is now defined over a larger
$p$-adic Lie extension. The rationale (which now seems irrational)
behind this thought is that our Selmer group is now a module over an
Iwasawa algebra of the group $\Zp^r\times(\Zp\rtimes\Zp)$ which has
a large ``commutative" component and, therefore, one might naively
hope that having large ``commutative" component may somehow force
the existence of nontrivial global annihilator for our Selmer group.
However, as our results (both algebraic and arithmetic) will show,
such an idea is not feasible in general.

\begin{ack}
    The work was first written up when the author was a Postdoctoral fellow
    at the GANITA Lab at the University of Toronto. He would like to acknowledge the
    hospitality and conducive working conditions provided by the GANITA
    Lab and the University of Toronto while this work was in progress.
    He would also like to thank the anonymous referee for his valuable
    comments, and for pointing out
    some mistakes in an earlier version of the paper.
            \end{ack}

\section{Algebraic Preliminaries} \label{Algebraic Preliminaries}

In this section, we establish some algebraic preliminaries and
notation. Throughout the paper, we will always work with left
modules over a ring. Let $\La$ be a (not neccessarily commutative)
Noetherian ring which has no zero divisors. Then it admits a skew
field of fractions $K(\La)$ which is flat over $\La$ (see
\cite[Chapters 6 and 10]{GW} or \cite[Chapter 4, \S 9 and \S
10]{Lam}). If $M$ is a finitely generated $\La$-module, we define
the $\La$-rank of $M$ to be
$$ \rank_{\La}M  = \dim_{K(\La)} K(\La)\ot_{\La}M. $$
Clearly, one has $\rank_{\La}M =0 $ if and only if
$K(\La)\ot_{\La}M=0$. We say that $M$ is a \textit{torsion}-module
if $\rank_{\La}M =0$. We shall record a simple lemma which is a
special case of \cite[Lemma 4.1]{LimFine}.

 \bl \label{torsion x} Let $x$ be a central element of $\La$ with
the property that $\Omega:=\La/x\La$ also has no zero divisors. Let
$M$ be a finitely generated $\La$-module. Then
  $$ \rank_{\Omega} M/xM = \rank_{\Omega}M[x] +
  \rank_{\La}M,
 $$ where $M[x]$ is the submodule of $M$ killed by $x$. \el

For a nonzero $\La$-module $M$, we define the global annihilator
ideal
\[ \Ann_{\La}(M) = \{ \la\in \La : \la m = 0~\mbox{for all}~m\in
M\}.\] Note that this is a two-sided ideal of $\La$. Indeed, for a
given $\la\in \Ann_{\La}(M)$, we have $\mu\la m = 0$ for every
$\mu\in\La$ and $m\in M$. On the other hand, since $\mu m \in M$, we
also have $\la\mu m = 0$. We will say that $M$ is a
\textit{faithful} $\La$-module if $\Ann_{\La}(M) = 0$.

Now if $x\in \La$, we denote $M[x]$ to be the set consisting of
elements of $M$ annihilated by $x$. If $x$ is not central, $M[x]$ is
at most an additive subgroup of $M$. However, if we assume further
that $x\La = \La x$, then it is easy to see that $M[x]$ is a
$\La$-submodule of $M$. Indeed, given $\la\in\La$ and $m\in M$, it
follows from the hypothesis $x\La = \La x$ that $\la'x=x\la$ for
some $\la'\in\La$, and therefore, one has $x\la m= \la' xm =0$.
Continuing to assume that $x\La = \La x$, one can also verify that
 \[xM = \{xm : m\in M\}  \]
is a $\La$-submodule of $M$, and that $M/xM$ is a $\La/x\La$-module.
We finally point out that $x\La$ is a two-sided ideal under the
condition that $x\La = \La x$.

The following lemma is a natural generalization of \cite[Lemma
4.5]{V03} and will be a crucial ingredient in proving our main
results in Sections 3 and 6. It will be of interest to have an
analogous statement for completely faithful modules (see Remark
\ref{control remark on complete faithful}) but we are not able to
establish such a statement at this point of writing.

 \bl \label{main lemma} Let $x$ be a nonzero element of $\La$ with
the property that $x\La = \La x$ and suppose that the ring
$\Omega:=\La/x\La$ has no zero divisors. Write $I= x\La $ $(=\La
x)$. Let $M$ be a finitely generated $\La$-module. Suppose that
$M[x]= 0$, and suppose that $\cap_{i\geq 1}I^i = 0$.

If $M/xM$ is a faithful $\Omega$-module, then $M$ is a faithful
$\La$-module. \el

\bpf
 We will prove the contrapositive statement. Suppose that
 $\Ann_{\La}(M)$ contains a nonzero element $\la$. Since $\cap_{i\geq 1}I^i =
 0$, we can find $n$ such that $\la\in I^n$ but $\la\notin
 I^{n+1}$. This in turn implies that $\la = x^n\la_0$ for some
 $\la_0\notin I$ (note that such a representation is possible
 by the assumption that $x\La = \La x$).
 But since $M[x]= 0$, we actually have $\la_0\in
 \Ann_{\La}(M)$. Since $\la_0\notin I$, the image of $\la_0$ in
 $\Omega$ is nonzero and lies in $\Ann_{\Omega}(M/xM)$.
\epf

We record two more lemmas. The first has an easy proof which is left
to the reader.

\bl \label{compare lemma pair} Suppose that we are given an exact
sequence
 \[0\lra B'\lra M' \lra M \lra B\lra 0\]
 of $\La$-modules, where $B'$ and $B$ are both annihilated by a
 nonzero central element
 $\la$ of $\La$. Then $M'$ is faithful over $\La$ if and only if $M$
 is faithful over $\La$. \el

\bl \label{compare lemma} Suppose that we are given an exact
sequence
 \[M' \lra M \lra M'' \lra 0\]
 of $\La$-modules, where $M'$ is nonzero and
 $M''$ is finite. Then if $M$ is faithful over $\La$, so is $M'$. \el

\bpf If the ring $\La$ has characteristic zero, one can probably
give a proof along the line of the proof of Lemma \ref{compare lemma
pair}. In fact, for much of the discussion in the paper, this case
will suffice. However, we thought that it may be of interest to give
a proof that works in general which we do now. We will prove the
contrapositive statement. Suppose that $\Ann_{\La}(M')$ contains a
nonzero element $\la$. Since $M''$ is finite, it follows that for
each $z \in M''$, there exists $n_z < m_z$ such that $\la^{n_z} z =
\la^{m_z}z$. This in turns implies that $\la^{n_z}(\la^{m_z-n_z}-
1)z = 0$. Set $n = 1 + \max_{z\in M''\setminus\{0\}} n_z$ and $m =
\prod_{z\in M''\setminus\{0\}} (m_z -n_z)$. Clearly, one has
$\la^n(\la^m -1)z = 0$ for every $z \in M''$. Also, since $n > 0$ by
our choice, we have that $\la^n(\la^m-1)$ lies in $\Ann_{\La}(M')$.
Therefore, $\la^n(\la^m-1)$ lies in $\Ann_{\La}(M)$. It remains to
show that $\la^n(\la^m-1)$ is a nonzero element of $\La$. Let $w$ be
a nonzero element of $M'$ (such an element exists by our hypothesis
that $M' \neq 0$). Then $(\la^m - 1)w = -w\neq 0$, and this in turn
implies that $\la^m - 1\neq 0$. Since $\La$ has no zero divisors and
$\la\neq 0$, it follows that $\la^n(\la^m-1)$ is also nonzero. This
completes the proof of the lemma. \epf

Let $\La$ be a Auslander regular ring (see \cite[Definition
3.3]{V02}) with no zero divisors. Let $M$ be a finitely generated
$\La$-module. Then $M$ is a torsion $\La$-module if and only if
$\Hom_{\La}(M,\La)=0$ (cf. \cite[Lemma 4.2]{LimFine}). If $M$ is a
torsion $\La$-module, we say that $M$ is a \textit{pseudo-null}
$\La$-module if $\Ext^1_{\La}(M,\La)=0$. Let $\mathcal{M}$ denote
the category of all finitely generated $\La$-modules, let
$\mathcal{C}$ denote the full subcategory of all pseudo-null modules
in $\mathcal{M}$ and let $q :\mathcal{M}\lra
\mathcal{M}/\mathcal{C}$ denote the quotient functor. For a finitely
generated $\La$-module $M$, we say that $M$ is \textit{completely
faithful} if  $Ann_{\La}(N) = 0$ for any $N \in\mathcal{M}$ such
that $q(N)$ is isomorphic to a non-zero subquotient of $q(M)$.

\bl \label{faithful ext} Let $\La$ be an Auslander regular ring with
no zero divisors. Then we have the following statements.

\begin{enumerate}
\item[$(a)$] If $M$ is completely faithful over $\La$, so is every non
pseudo-null subquotient of $M$.

\item[$(b)$] An extension of completely faithful $\La$-modules is also
completely faithful.
\end{enumerate}
\el

\bpf
 This is straightforward from the definition.
\epf

\section{Completely faithful modules over completed group algebras}

In this section, we will prove our main theorem which is an
extension of \cite[Theorem 6.3]{V03}. As before, $p$ will denote a
fixed odd prime. Let $G$ be a compact pro-$p$ $p$-adic Lie group
without $p$-torsion. It is well known that $\Zp\llbracket
G\rrbracket$ is an Auslander regular ring (cf. \cite[Theorems
3.26]{V02}). Furthermore, the ring $\Zp\ps{G}$ has no zero divisors
(cf.\ \cite{Neu}), and therefore, as seen in the previous section,
there is a well-defined notion of $\Zp\ps{G}$-rank and torsion
$\Zp\ps{G}$-module. We record the following well-known and important
result of Venjakob (cf. \cite[Example 2.3 and Proposition
5.4]{V03}).

\bt[Venjakob] \label{pseudo-null torsion} Suppose that $H$ is a
closed normal subgroup of $G$ with $G/H \cong \Zp$. Let $M$ be a
compact $\Zp\ps{G}$-module which is finitely generated over
$\Zp\ps{H}$. Then $M$ is a pseudo-null $\Zp\ps{G}$-module if and
only if $M$ is a torsion $\Zp\ps{H}$-module. \et

We record another useful lemma whose proof is left to the reader (or
see \cite[Lemma 4.5]{LimFine}).

\bl \label{relative rank} Let $H$ be a compact pro-$p$ $p$-adic Lie
group without $p$-torsion. Let $N$ be a closed normal subgroup of
$H$ such that $N\cong \Zp$ and such that $H/N$ is also a compact
pro-$p$ $p$-adic Lie group without $p$-torsion. Let $M$ be a
finitely generated $\Zp\ps{H}$-module. Then $H_i(N,M)$ is finitely
generated over $\Zp\ps{H/N}$ for each $i$ and $H_i(N,M) =0$ for
$i\geq 2$. Furthermore, we have an equality
\[\rank_{\Zp\ps{H}}M = \rank_{\Zp\ps{H/N}}M_N - \rank_{\Zp\ps{H/N}}H_1(N,M). \] \el

Before continuing our discussion, we introduce the following
hypothesis on our group $G$.

\medskip \noindent
$\mathbf{(NH)}$ : The group $G$ contains two closed normal subgroups
$N$ and $H$ which satisfy the following two properties.
  \begin{enumerate}
\item[$(i)$] $N\subseteq H$, $G/H\cong \Zp$ and
 $G/N$ is a non-abelian group isomorphic to $\Zp\rtimes\Zp$.
\item[$(ii)$] There is a finite family of closed normal subgroups $N_i$
$(0\leq i\leq r)$ of $G$ such that $1=N_0\subseteq N_1 \subseteq
\cdots\subseteq N_r =N$ and such that $N_i/N_{i-1}\cong \Zp$ for
$1\leq i\leq r$.
\end{enumerate}

We can now state and prove the main theorem of this section which
generalizes \cite[Theorem 6.3]{V03}.

\bt \label{faithful modules}
  Suppose that $G$ satisfies  $\mathbf{(NH)}$.
 Let $M$ be a $\Zp\ps{G}$-module which
is finitely generated over $\Zp\ps{H}$ and has positive
$\Zp\ps{H}$-rank. Then $M$ is a completely faithful
$\Zp\ps{G}$-module. \et

\bpf
 Since $M$ is finitely generated over $\Zp\ps{H}$, it follows
 from the result of Venjakob mentioned above that every subquotient of $M$
 which is not pseudo-null has positive $\Zp\ps{H}$-rank. Therefore,
 it suffices to show that every $\Zp\ps{G}$-module that
is finitely generated over $\Zp\ps{H}$ with positive
$\Zp\ps{H}$-rank is faithful over $\Zp\ps{G}$. We will proceed by
induction on $r$. When $r=0$, this is precisely \cite[Corollary
4.3]{V03}. Now suppose that $r\geq 1$ and suppose that $M$ is
finitely generated over $\Zp\ps{H}$ with positive $\Zp\ps{H}$-rank.
Then choose a topological generator $\ga_1$ of $N_1$. The ideal
generated by $\ga_1 -1$ is precisely the augmentation kernel
$I_{N_1}$ of the canonical quotient map $\Zp\ps{G}\twoheadrightarrow
\Zp\ps{G/N_1}$, and one has $I_{N_1} =(\ga_1 -1)\Zp\ps{G}=
\Zp\ps{G}(\ga_1 -1)$. Therefore, $M[\ga_1 -1]$ is a
$\Zp\ps{G}$-submodule of $M$. Suppose for now that $M[\ga_1-1] =0$.
Then it follows from Lemma \ref{relative rank} that
 \[ \rank_{\Zp\ps{H/N_1}} M_{N_1} =  \rank_{\Zp\ps{G}} M  +
 \rank_{\Zp\ps{H/N_1}} H_1(N_1,M)
  =  \rank_{\Zp\ps{H}} M >0. \]
  Here the second equality follows from the facts that $H_1(N_1,M) =
  M[\ga_1-1]$ and that $M[\ga_1-1] =0$.
Hence  $M_{N_1}$ is a $\Zp\ps{G/N_1}$-module which is finitely
generated over $\Zp\ps{H/N_1}$ with positive $\Zp\ps{H/N_1}$-rank.
By our induction hypothesis, we have that $M_{N_1}$ is (completely)
faithful over $\Zp\ps{G/N_1}$.  Note that $I_{N_1}$ is closed in
$\Zp\ps{G}$ and so $\cap_{i\geq 1}I_{N_1}^i = 0$. Also, since we are
assuming that $M[\ga_1-1] =0$, we may apply Lemma \ref{main lemma}
to conclude that $M$ is faithful over $\Zp\ps{G}$.

It remains to consider the situation when $M[\ga_1-1] \neq 0$. Since
$(\ga_1 -1)^i\Zp\ps{G}= \Zp\ps{G}(\ga_1 -1)^i$ for all $i \geq 1$,
it follows that $M[(\ga_1-1)^i]$ is a $\Zp\ps{G}$-submodule of $M$
for every $i \geq 1$. As $M$ is a Noetherian $\Zp\ps{G}$-module, the
chain
\[ M[\ga_1-1] \subseteq M[(\ga_1-1)^2] \subseteq
\cdots \subseteq M[(\ga_1-1)^i]  \subseteq \cdots\]
 terminates at a finite level which in turn implies that
 \[ M' : = \cup_{i\geq 1}M[(\ga_1-1)^i] = M[(\ga_1-1)^n] \]
 for some $n$. Set $M'' =
 M/M'$. Since $M' = M[(\ga_1 -1)^n]$ is a torsion $\Zp\ps{H}$-module,
it follows that $M''$ has positive $\Zp\ps{H}$-rank. On the other
hand, one clearly has $\Ann_{\Zp\ps{G}}(M) \subseteq
\Ann_{\Zp\ps{G}}(M'')$, and therefore, we are reduced to showing
that $M''$ is faithful over $\Zp\ps{G}$. By our construction of
$M''$, we have that $M''[\ga_1-1]= 0$. Hence we may apply the
argument in the previous paragraph to obtain the faithfulness of
$M''$. The proof of the theorem is now complete. \epf

\br
  One can of course prove analogous result as in Theorem \ref{faithful
  modules} replacing the coefficient ring $\Zp$ by $\mathbb{F}_p$.
  This is achieved by combining
  the argument in Theorem \ref{faithful
  modules} with \cite[Proposition 4.2(i)]{V03}.
  In fact, it is not difficult (though not immediate) to see that one can also prove
  analogous results as in the paper \cite{V03} and
  Theorem \ref{faithful modules} for a ring of integer
  of a finite extension of $\Qp$ and its residue field.
\er

We should mention here that the conclusion in Theorem \ref{faithful
modules} does not hold for a general $G$. We give a class of
counterexamples. For the remainder of this paragraph, we assume that
$G =H\times G/H$. Let $M=\Zp\ps{H}$ be the $\Zp\ps{G}$-module where
the action of $H$ is the natural one, and the action of $G/H$ is the
trivial one. Denote $\ga$ to be a topological generator of $G/H$.
Clearly, $M$ is clearly annihilated by $\ga-1$, and so it is not
faithful over $\Zp\ps{G}$.

Despite the counterexamples, it is still of interest to ask if there
exists other torsionfree pro-$p$ $p$-adic Lie group $G$ where a
similar conclusion in Theorem \ref{faithful modules} holds.
Alternatively, one may ask if $G$ is a torsionfree pro-$p$ $p$-adic
Lie group with a normal subgroup $H$ such that $G/H\cong \Zp$, and
has the property that every finitely generated $\Zp\ps{H}$-module of
positive $\Zp\ps{H}$-rank is completely faithful over $\Zp\ps{G}$,
what can one say about the structure of $G$? The author does not
have an answer to these questions at this point of writing.

Finally, we end the section mentioning how the results in this
section can be extended to a certain class of finitely generated
torsion $\Zp\ps{G}$ which was first introduced in \cite{CFKSV}. In
particular, this class of modules is a source of examples of
\textit{faithful modules that are not completely faithful}. As
before, $G$ is a compact pro-$p$ $p$-adic Lie group without
$p$-torsion and $H$ is a closed normal subgroup of $G$ such that
$G/H\cong \Zp$. For a finitely generated $\Zp\ps{G}$-module $M$, we
say that $M$ \textit{belongs to} $\M_H(G)$ if $M/M(p)$ is finitely
generated over $\Zp\ps{H}$. Here $M(p)$ is the submodule of $M$
consisting of elements of $M$ annihilated by a power of $p$. Note
that a finitely generated $\Zp\ps{G}$-module belonging to $\M_H(G)$
is necessarily a torsion $\Zp\ps{G}$-module.

\bc \label{faithful modules MHG}
  Suppose that $G$ satisfies $\mathbf{(NH)}$.
 Let $M$ be a $\Zp\ps{G}$-module which
belongs to $\M_H(G)$ and has the property that $M/M(p)$ has a
positive $\Zp\ps{H}$-rank. Then $M$ is a faithful
$\Zp\ps{G}$-module.

Furthermore, $M$ is a completely faithful $\Zp\ps{G}$-module if and
only if $M(p)$ is a pseudo-null $\Zp\ps{G}$-module. \ec

\bpf
 Since $M$ is finitely generated over
$\Zp\ps{G}$, the module $M(p)$ is annihilated by a power of $p$. The
first assertion is now an immediate consequence from Lemma
\ref{compare lemma pair} and Theorem \ref{faithful modules}.

To prove the second assertion, we first recall that
$q:\mathcal{M}\lra \mathcal{M}/\mathcal{C}$ is the quotient functor,
where $\mathcal{M}$ denotes the category of all finitely generated
$\La$-modules and $\mathcal{C}$ denotes the full subcategory of all
pseudo-null modules in $\mathcal{M}$. Now suppose that $M(p)$ is a
pseudo-null $\Zp\ps{G}$-module, then we have $q(M) = q
\big(M/M(p)\big)$. Since $M/M(p)$ is completely faithful by Theorem
\ref{faithful modules}, it follows that $M$ is also completely
faithful. On the other hand, if $M(p)$ is not a pseudo-null
$\Zp\ps{G}$-module and is annihilated by some power of $p$, then $M$
contains a submodule which is not pseudo-null and not faithful, and
therefore, is not completely faithful. \epf

\br
 Note that $M(p)$ is a pseudo-null $\Zp\ps{G}$-module if and only if its
$\mu_G$-invariant (see \cite[Definition 3.32]{V02} for definition)
vanishes (cf. \cite[Remark 3.33]{V02}). \er

\section{Completely faithful Selmer groups} \label{completely
faithful Sel section}

Let $F$ be a number field. Fix once and for all an algebraic closure
$\bar{F}$ of $F$. Therefore, an algebraic (possibly infinite)
extension of $F$ will mean an subfield of $\bar{F}$ which contains
$F$. Let $E$ be an elliptic curve over $F$. Assume that for every
prime $v$ of $F$ above $p$, our elliptic curve $E$ has either good
ordinary reduction or multiplicative reduction at $v$.

 Let $v$ be a prime of $F$. For every finite extension $L$ of
$F$, we define
 \[ J_v(E/L) = \bigoplus_{w|v}H^1(L_w, E)_{p^{\infty}},\]
where $w$ runs over the (finite) set of primes of $L$ above $v$. If
$\mathcal{L}$ is an infinite extension of $F$, we define
\[ J_v(E/\mathcal{L}) = \ilim_L J_v(E/L),\]
where the direct limit is taken over all finite extensions $L$ of
$F$ contained in $\mathcal{L}$. For any algebraic (possibly
infinite) extension $\mathcal{L}$ of $F$, the Selmer group of $E$
over $\mathcal{L}$ is defined to be
\[ S(E/\mathcal{L}) = \ker\Big(H^1(\mathcal{L}, E_{p^{\infty}})
\lra \bigoplus_{v} J_v(E/\mathcal{L}) \Big), \] where $v$ runs
through all the primes of $F$.

We say that $F_{\infty}$ is an \textit{admissible $p$-adic Lie
extension} of $F$ if (i) $\Gal(F_{\infty}/F)$ is a compact $p$-adic
Lie group, (ii) $F_{\infty}$ contains the cyclotomic $\Zp$-extension
$F^{\cyc}$ of $F$ and (iii) $F_{\infty}$ is unramified outside a
finite set of primes of $F$. Furthermore, an admissible $p$-adic Lie
extension $F_{\infty}$ of $F$ will be said to be \textit{strongly
admissible} if $\Gal(F_{\infty}/F)$ is a compact pro-$p$ $p$-adic
Lie group without $p$-torsion. Write $G = \Gal(F_{\infty}/F)$, $H =
\Gal(F_{\infty}/F^{\cyc})$ and $\Ga =\Gal(F^{\cyc}/F)$. Let $S$ be a
finite set of primes of $F$ which contains the primes above $p$, the
infinite primes, the primes at which $E$ has bad reduction and the
primes that are ramified in $F_{\infty}/F$. Denote $F_S$ to be the
maximal algebraic extension of $F$ unramified outside $S$. For each
algebraic (possibly infinite) extension $\mathcal{L}$ of $F$
contained in $F_S$, we write $G_S(\mathcal{L}) =
\Gal(F_S/\mathcal{L})$. The following alternative equivalent
description of the Selmer group of $E$ over $F_{\infty}$
\[ S(E/F_{\infty}) = \ker\Big(H^1(G_S(F_{\infty}), E_{p^{\infty}})
\stackrel{\lambda_{E/F_{\infty}}}{\lra} \bigoplus_{v\in S}
J_v(E/F_{\infty}) \Big)\] is well-known (for instance, see
\cite[Lemma 2.2]{CS12}). We will denote $X(E/F_{\infty})$ to be the
Pontryagin dual of $S(E/F_{\infty})$. The following is then an
immediate consequence of Theorem \ref{faithful modules}.

\bt \label{faithful Selmer groups}
 Let $E$ be an elliptic curve over $F$
which has either good ordinary reduction or multiplicative reduction
at every prime of $F$ above $p$.
 Let $F_{\infty}$ be a strongly admissible $p$-adic Lie extension of $F$ with
 $G =\Gal(F_{\infty}/F)$. Suppose that $G$ satisfies $\mathbf{(NH)}$.
  If $X(E/F_{\infty})$ is finitely generated over
 $\Zp\ps{H}$ with positive $\Zp\ps{H}$-rank, then $X(E/F_{\infty})$
 is completely faithful over $\Zp\ps{G}$.
\et

We record the following corollary of the theorem which is useful in
obtaining examples of completely faithful Selmer groups.

\bc \label{faithful corollary}
 Let $F_{\infty}$ be a strongly admissible $p$-adic Lie extension of $F$ with
 $G=\Gal(F_{\infty}/F)$. Suppose that $G$ satisfies $\mathbf{(NH)}$.
 Assume that $X(E/F^{\cyc})$ is finitely generated over $\Zp$.
 Furthermore, suppose that either of the following conditions is
 satisfied.
 \begin{enumerate}
 \item[$(a)$]
  $X(E/F^{\cyc})$ has positive $\Zp$-rank.
 \item[$(b)$]
  The field $F$ is not totally real, the elliptic curve $E$ has good
  ordinary reduction at every prime of $F$
  above $p$ and $X(E/F_{\infty})\neq 0$.
 \end{enumerate}
  Then $X(E/F_{\infty})$ is completely faithful over $\Zp\ps{G}$. \ec

\bpf
 By \cite[Proposition 5.6]{CFKSV}
or \cite[Theorem 2.1]{CS12}, if $X(E/F^{\cyc})$ is finitely
generated over $\Zp$, then $X(E/F_{\infty})$ is finitely generated
over $\Zp\ps{H}$. It suffices to show that $X(E/F_{\infty})$ has
 positive $\Zp\ps{H}$-rank under the assumption of either
 conditions. If (a) holds, a standard argument in the spirit of
 \cite{HM} will allow us establish the positivity of
 $\Zp\ps{H}$-rank (alternatively, one can make use of \cite[Theorem 5.4]{HS} directly).
 Now if (b) holds, we may apply the main result in \cite{Mat}
 to conclude that $X(E/F_{\infty})$ has
 positive $\Zp\ps{H}$-rank (for instance, see \cite[Lemma 5.8]{LimMHG}).
\epf

\br Of course, one can apply Theorem \ref{faithful modules} to
obtain completely faithful Selmer groups of $p$-ordinary modular
forms, or even more general $p$-adic representations. \er

We now discuss the complete faithfulness of Selmer groups of Hida
deformations. Let $E$ be an elliptic curve over $\Q$ with ordinary
reduction at $p$ and assume that $E[p]$ is an absolutely irreducible
$\Gal(\bar{\Q}/\Q)$-representation. By Hida theory (for instance,
see \cite{Hi1, Hi2}), there exists a commutative complete Noetherian
local domain $R$ which is flat over the power series ring
$\Zp\ps{X}$ in one variable, and a free $R$-module $T$ of rank 2
with $T/P \cong T_pE$ for some prime ideal $P$ of $R$. \textit{We
will further assume that $R = \Zp\ps{X}$ in all our discussion}. For
more detailed description of fundamental and important arithmetic
properties of the Hida deformations, we refer readers to \cite{CS12,
Hi1, Hi2, Jh, SS}. We will just mention two properties of $T$ which
we require to define an appropriate Selmer group of the Hida
deformation. The first is that $T$ is unramified outside the set
$S$, where $S$ is any finite set of primes of $F$ which contains the
primes above $p$, the infinite primes, the primes at which $E$ has
bad reduction and the primes that are ramified in $F_{\infty}/F$.
The second property we will mention is that there exists an
$R$-submodule $T^+$ of $T$ which is invariant under the action of
$\Gal(\bar{\Q}_p/\Qp)$ and such that both $T^+$ and $T/T^+$ are free
$R$-modules of rank one.

Set $A = T\ot_R\Hom_{\cts}(R,\Qp/\Zp)$ and $A^+ =
T^+\ot_R\Hom_{\cts}(R,\Qp/\Zp)$. We note that one has
$E_{p^{\infty}} = A[P]$. Then following \cite[Section 4]{CS12} or
\cite[Section 6]{SS}, we define the Selmer group of the Hida
deformation over an admissible $p$-adic Lie extension $F_{\infty}$
of $\Q$ by

\[  S(A/F_{\infty}) = \ker\Big(H^1(G_S(F_{\infty}), A)
\lra \bigoplus_{v\in S} J_v(A, F_{\infty}) \Big), \]
 where
 \[ J_v(A, F_{\infty}) = \begin{cases}
  \prod_{w|v}H^1(F_{\infty,w}, A/A^+),& \mbox{if } v\mbox{ divides }p, \\
      \prod_{w|v}H^1(F_{\infty,w}, A), &
      \mbox{if } v\mbox{ does not divides }p. \end{cases}
 \]
 We will denote by $X(A/F_{\infty})$ the Pontryagin dual of this
Selmer group. We will consider this dual Selmer group as a (compact)
$\Gal(F_{\infty}/F)$-module for some finite extension $F$ of $\Q$ in
$F_{\infty}$, where $F_{\infty}$ is a strongly admissible $p$-adic
Lie extension of $F$.

\bt \label{faithful Selmer groups for Hida}
 Let $F_{\infty}$ be a strongly admissible $p$-adic Lie extension of $F$ with
 Galois group $G$. Suppose that $G$ satisfies $\mathbf{(NH)}$.
 If $X(A/F_{\infty})$ is finitely generated over
 $R\ps{H}$ with positive $R\ps{H}$-rank, then $X(A/F_{\infty})$
 is completely faithful over $R\ps{G}$. In particular, if $X(E/F_{\infty})$
 is finitely generated over
 $\Zp\ps{H}$ with positive $\Zp\ps{H}$-rank, then $X(A/F_{\infty})$
 is completely faithful over $R\ps{G}$.
\et

\bpf
 By identifying $R\ps{G} \cong
 \Zp\ps{\Zp\times G}$, the first part of the theorem is
immediate from Theorem \ref{faithful modules}. For the second part,
it suffices to show that $X(A/F_{\infty})$ is finitely generated
over $R\ps{H}$ with positive $R\ps{H}$-rank. By \cite[Theorem
4.2]{CS12}, the map
 \[ X(A/F_{\infty})/P \lra X(E/F_{\infty}) \]
has cokernel which is finitely generated over $\Zp$. Since
$X(E/F_{\infty})$ has positive $\Zp\ps{H}$-rank, this in turn
implies that $X(A/F_{\infty})/P$ has positive $\Zp\ps{H}$-rank. By
an application of an argument in the spirit to that in \cite[Theorem
7.4]{SS}, one can show that $S(A/F_{\infty})/P =0$. Equivalently,
this is the same as saying that $X(A/F_{\infty})[P] =0$. Therefore,
we may apply Lemma \ref{torsion x} to conclude that $R\ps{H}$-rank
of $X(A/F_{\infty})$ is the same as the $\Zp\ps{H}$-rank of
$X(E/F_{\infty})$ and, in particularly, is positive as required.
\epf

\br One can also obtain completely faithful Selmer groups of
$p$-adic representations defined over coefficient rings $\Zp\ps{X_1,
X_2, ..., X_n}$ over strongly admissible $p$-adic Lie extensions of
the form considered in this section. \er

We now give some examples to illustrate the results in this section.

\medskip
(a) Let $E$ be the elliptic curve $11a2$ of Cremona's table which is
given by
 \[ y^2 + y = x^3 -x. \]
Take $p = 5$, $F=\Q(\mu_5)$ and $L_{\infty} = \Q(\mu_{5^{\infty}},
11^{5^{-\infty}})$. By \cite[Theorem 6.2]{HV}, $X(E/L_{\infty})$ is
a free $\Z_5\ps{\Gal(L_{\infty}/F^{\cyc})}$-module of rank four. Let
$F_{\infty}$ be a strongly admissible $5$-adic Lie extension of $F$
that contains $L_{\infty}$ and that the group
$N=\Gal(F_{\infty}/L_{\infty})$ satisfies the conditions in Theorem
\ref{faithful Selmer groups}. A similar argument as in Corollary
\ref{faithful corollary} shows that $X(E/F_{\infty})$ is a finitely
generated $\Z_5\ps{H}$-module of positive rank. Then we may apply
Theorem \ref{faithful Selmer groups} to conclude that
$X(E/F_{\infty})$ is completely faithful over $\Z_5\ps{G}$. Now if
$E'$ is either $11a1$ or $11a3$, then it follows from Proposition
\ref{faithful isogeny} that $X(E'/F_{\infty})$ is faithful over
$\Z_5\ps{G}$. We claim that in either cases, $X(E'/F_{\infty})$ is
not completely faithful. To see this, it suffices, by Proposition
\ref{faithful modules MHG}, to show that $X(E'/F_{\infty})(5)$ is
not pseudo-null, or equivalently, the $\mu_G$-invariant of
$X(E'/F_{\infty})$ is positive. By an application of \cite[Lemma
5.6]{CFKSV}, we have that $X(E'/F_{\infty})$ belongs to $\M_H(G)$.
This in turn allows us to apply \cite[Theorem 3.1]{LimMHG} to
conclude that $\mu_G(X(E'/F_{\infty})) = \mu_{\Ga}(X(E'/F_{\cyc}))$.
But this latter quantity is well-known to be nonzero (cf.\
\cite[Theorem 5.28]{CS00}) and hence our claim is established.

\medskip (b) The next example is taken from \cite{Jh}. Let $E$ be the
elliptic curve $79a1$ of Cremona's table which is given by
\[ y^2 + xy + y = x^3 + x^2 -2x.\]
 Take $p = 3$ and $F=\Q(\mu_3)$. As
noted in \cite{Jh}, $X(E/F^{\cyc})$ is isomorphic to $\Z_3$.  Let
$F_{\infty}$ be a strongly admissible $3$-adic Lie extension of $F$
that satisfies the conditions in Theorem \ref{faithful Selmer
groups}. Write $G = \Gal(F_{\infty}/F)$. By Corollary \ref{faithful
corollary}, $X(E/F_{\infty})$ is a completely faithful
$\Z_3\ps{G}$-module. Let $A$ be the Galois module obtained from the
Hida family associated to $E$ as above. Therefore, one may apply
Theorem \ref{faithful Selmer groups for Hida} to conclude that
$X(A/F_{\infty})$ is a completely faithful $R\ps{G}$-module.
Examples of strongly admissible $3$-adic extensions of $F$ that one
can take are:
\[\Q(\mu_{3^{\infty}}, 2^{3^{-\infty}}),\quad
\Q(\mu_{3^{\infty}}, 2^{3^{-\infty}}, 3^{3^{-\infty}}),
\quad\Q(\mu_{3^{\infty}}, 3^{3^{-\infty}}, 5^{3^{-\infty}}),\quad
\Q(\mu_{3^{\infty}}, 2^{3^{-\infty}}, 3^{3^{-\infty}},
5^{3^{-\infty}}), \quad... ~\mbox{etc}.\]

\section{Isogeny invariance of faithful Selmer
groups} \label{isogeny section}

In this short section, we will show that the property of
faithfulness is an isogeny invariant. Namely, we prove the following
statement.

\bp \label{faithful isogeny}
 Let $E_1$  and $E_2$ be two elliptic curves over $F$
with either good ordinary reduction or multiplicative reduction at
every prime of $F$ above $p$ which are isogenous to each other. Let
$F_{\infty}$ be a strongly admissible noncommutative $p$-adic Lie
extension of $F$ with $G= \Gal(F_{\infty}/F)$. Assume that both
$X(E_1/F_{\infty})$ and $X(E_2/F_{\infty})$ are torsion over
$\Zp\ps{G}$, and that the localization maps
$\lambda_{E_1/F_{\infty}}$ and $\lambda_{E_2/F_{\infty}}$ are
surjective. Then $X(E_1/F_{\infty})$ is a faithful
$\Zp\ps{G}$-module if and only if $X(E_2/F_{\infty})$ is a faithful
$\Zp\ps{G}$-module. \ep

\bpf Let $\varphi: E_1\lra E_2$ be an isogeny defined over $F$.
 By a standard argument to that in the proof of \cite[Theorem 5.1]{HV} or
 \cite[Theorem 3.1]{Ho}, we can show that $\varphi$ induces a
 $\Zp\ps{G}$-homomorphism
 \[ X(E_2/F_{\infty})\lra X(E_1/F_{\infty}),\]
 whose kernel and cokernel are killed by $p^n$ for some large enough
 $n$. The required conclusion is now immediate from an application
 of Lemma \ref{compare lemma pair}.
\epf

On the other hand, completely faithfulness is not an isogeny
invariant. As seen in the previous section, the dual Selmer group of
11a2 is completely faithful but the the dual Selmer group of 11a1
and 11a3 are not.

\section{Control theorems for faithfulness of Selmer
groups} \label{control section}

In this section, we will prove two control theorems on faithfulness
of Selmer groups which can be applied to a $p$-adic Lie extension
whose Galois group does not satisfy $\mathbf{(NH)}$. In particular,
our result will show that one cannot obtain nonfaithful Selmer
groups from a faithful Selmer group by adjoining $\Zp^r$-extension
or moving into the Hida deformation in general. We retain the
notation of the Section \ref{completely faithful Sel section}.
Recall that $S$ is a finite set of primes of $F$ which contains the
primes above $p$, the infinite primes, the primes at which $E$ has
bad reduction and the primes that ramify in $F_{\infty}/F$. We now
record a lemma.

\bl \label{faithful control lemma}
 Let $E$ be an elliptic curve over $F$
with either good ordinary reduction or multiplicative reduction at
every prime of $F$ above $p$. Let $F_{\infty}$ be a strongly
admissible $p$-adic Lie extension of $F$ with $G=
\Gal(F_{\infty}/F)$. Let $N$ be a closed normal subgroup of $G$ such
that $N\cong \Zp$ and that $G/N$ is a pro-$p$ $p$-adic Lie group
without $p$-torsion. Set $L_{\infty} = F_{\infty}^N$. Consider the
following statements.
\begin{enumerate}
 \item[$(i)$] $X(E/L_{\infty})$ is a torsion $\Zp\ps{G/N}$-module.
 \item[$(ii)$] $E_{p^{\infty}}$ is not rational over $L_{\infty}$.
 \item[$(iii)$] For every $v\in S$, the decomposition group of $G$
  at $v$ has dimension $\geq 2$.
\end{enumerate}
If $(i)$ holds, then $X(E/F_{\infty})$ is a torsion
$\Zp\ps{G}$-module. If $(i)$ and $(ii)$ hold, or $(i)$ and $(iii)$
hold, then we have that $H_1(N, X(E/F_{\infty})) = 0$. \el

\bpf
 By an argument similar to that in
\cite[Lemma 2.4]{CS12}, we have that the dual of the restriction map
  \[  X(E/F_{\infty})_N \lra X(E/L_{\infty}) \]
has kernel and cokernel which are finitely generated over
$\Zp\ps{H/N}$. Therefore, if (i) holds, it then follows that
$X(E/F_{\infty})_N$ is a torsion $\Zp\ps{G/N}$-module. By Lemma
\ref{relative rank}, this in turns implies that $X(E/F_{\infty})$ is
a torsion $\Zp\ps{G}$-module.

Now if we suppose that either (i) and (ii) hold, or (i) and (iii)
hold. Then by \cite[Proposition 3.3]{LimMHG}, we have a short exact
sequence
\[ 0 \lra S(E/L_{\infty}) \lra H^1(G_S(L_{\infty}), E_{p^{\infty}})
\stackrel{\lambda_{E/L_{\infty}}}{\lra} \bigoplus_{v\in S}
J_v(E/L_{\infty}) \lra 0 \]
 and that $H^2(G_S(L_{\infty}), E_{p^{\infty}})=0$.
By \cite[Lemma 7.1, Proposition 7.2]{LimFine}, the latter in turn
implies that $H^2(G_S(F_{\infty}), E_{p^{\infty}})=0$. Now we may
apply a similar argument to that in \cite[Theorem 7.4]{SS} (or
\cite[Lemma 4.4, Proposition 4.5]{LimMHG}) to obtain the conclusion
that $H_1(N, X(E/F_{\infty})) = 0$. \epf

\br
 We can obtain $H_1(N, X(E/F_{\infty})) = 0$ without assumptions
 (ii) and (iii), provided one assumes that $X(E/L_{\infty})$ belongs to
 $\M_{H/N}(G/N)$. To see this, we first note that since
 $X(E/L_{\infty})$ belongs to
 $\M_{H/N}(G/N)$, it follows from \cite[Proposition 2.5]{CS12} that,
 for every finite extension $K$
of $F$ in $L_{\infty}$, $X(E/K^{\cyc})$ is $\Zp\ps{\Ga_K}$-torsion,
where $\Ga_K=\Gal(K^{\cyc}/K)$. By \cite[Corollary 3.4]{LimMHG},
this in turn implies that we have a short exact sequence
\[ 0 \lra S(E/L_{\infty}) \lra H^1(G_S(L_{\infty}), E_{p^{\infty}})
\stackrel{\lambda_{E/L_{\infty}}}{\lra} \bigoplus_{v\in S}
J_v(E/L_{\infty}) \lra 0 \]
 and $H^2(G_S(L_{\infty}), E_{p^{\infty}})=0$. Now we may
proceed as in the argument to that in Lemma \ref{faithful control
lemma} to obtain the conclusion that $H_1(N, X(E/F_{\infty})) = 0$.
\er

We can now prove our first control theorem which concerns extensions
of admissible $p$-adic Lie extensions.

\bp \label{faithful control}
 Let $E$ be an elliptic curve over $F$
with either good ordinary reduction or multiplicative reduction at
every prime of $F$ above $p$.
 Let $F_{\infty}$ be a strongly admissible $p$-adic Lie extension of $F$ with
 $G= \Gal(F_{\infty}/F)$.
  Suppose that the following statements hold.
  \begin{enumerate}
\item[$(i)$] $N$ is a closed normal subgroup of $G$ which is contained
in $H$, and there is a finite family of closed normal subgroups
$N_i$ $(0\leq i\leq r)$ of $G$ such that $1=N_0\subseteq N_1
\subseteq \cdots\subseteq N_r =N$ and such that $N_i/N_{i-1}\cong
\Zp$ for $1\leq i\leq r$.
\item[$(ii)$] $G/N$ is a non-abelian pro-$p$ $p$-adic Lie group without
$p$-torsion. $($In particular, the dimension of the $p$-adic Lie
group $G/N$ is necessarily $\geq 2$.$)$
\item[$(iii)$] Set $L_{\infty} : = F_{\infty}^N$. Suppose that
$X(E/L_{\infty})$ is torsion over $\Zp\ps{G/N}$.
\item[$(iv)$] Suppose that
either $(a)$ or $(b)$ holds.
 \begin{enumerate}
 \item[$(a)$] $E_{p^{\infty}}$ is not rational over $F_{\infty}^{N_1}$ and
 $X(E/L_{\infty})$ is a faithful $\Zp\ps{G/N}$-module.
 \item[$(b)$] For every $v\in S$, the decomposition group of $\Gal(L_{\infty}/F)$
  at $v$ has dimension $\geq 2$,
  $X(E/L_{\infty})$ is a completely faithful $\Zp\ps{G/N}$-module
  and $r=1$.
\end{enumerate}
\end{enumerate}
 Then $X(E/F_{\infty})$ is a faithful $\Zp\ps{G}$-module.
\ep

\bpf
 We first note that by an iterative application of Lemma \ref{faithful control
lemma}, it follows from condition (i) and (iii) that
$X(E/F_{\infty})$ is torsion over $\Zp\ps{G}$. (Alternatively, one
may apply an argument similar to that in \cite[Theorem 2.3]{HO},
noting that $N$ is a solvable uniform pro-$p$ group.) Now suppose
that condition (iv)(a) holds. By an induction on $r$, it suffices to
prove the proposition in this case assuming $r=1$. By another
application of Lemma \ref{faithful control lemma}, one has that
$H_1(N, X(A/F_{\infty})) =0$. This in turn implies that
$X(E/F_{\infty})[\ga_N-1] = 0$, where $\ga_N$ is a topological
generator of $N$. As observed in the proof of Theorem \ref{faithful
modules}, $I_{N}$ is closed in $\Zp\ps{G}$ and so $\cap_{i\geq
1}I_{N_1}^i = 0$. Therefore, by Lemma \ref{main lemma}, we are
reduced to proving that $X(E/F_{\infty})_N$ is a faithful
$\Zp\ps{G/N}$-module. Now, applying an argument similar to that in
\cite[Lemma 2.4]{CS12}, we have that the dual of the cokernel of the
map
  \[ \al: X(E/F_{\infty})_N \lra X(E/L_{\infty}) \]
is contained in $H^1(N, E(F_{\infty})_{p^{\infty}})$. By the first
assumption in condition (iv)(a) and \cite[Proposition 10]{Z04}, we
have that $H^0(N, E(F_{\infty})_{p^{\infty}})=
E(L_{\infty})_{p^{\infty}}$ is finite. On the other hand, it follows
from Lemma \ref{relative rank} that
  \[ \corank_{\Zp}
E(L_{\infty})_{p^{\infty}} =\corank_{\Zp} H^1(N,
E(F_{\infty})_{p^{\infty}}) +\corank_{\Zp\ps{N}}
E(F_{\infty})_{p^{\infty}}.\]
 Therefore, it follows that $H^1(N,
E(F_{\infty})_{p^{\infty}})$, and hence the cokernel of $\al$, is
finite. We may now combine Lemma \ref{compare lemma} with the second
assumption of condition (iv)(a) to conclude that $X(E/F_{\infty})_N$
is a faithful $\Zp\ps{G/N}$-module.

We now consider the case when condition (iv)(b) holds. As above, it
suffices to show that $X(E/F_{\infty})_N$ is a faithful
$\Zp\ps{G/N}$-module. Again, by the argument of \cite[Lemma
2.4]{CS12}, one can show that the map
  \[ \al: X(E/F_{\infty})_N \lra X(E/L_{\infty}) \]
has kernel which is cofinitely generated over $\Zp\ps{H/N}$, and
cokernel which is cofinitely generated over $\Zp$. In particular, by
condition (ii) and Theorem \ref{pseudo-null torsion}, the dual of
the cokernel of $\al$ is psuedo-null over $\Zp\ps{G/N}$.
Furthermore, in view of the first assumption of condition (iv)(b),
one can apply a similar argument in the spirit of the proof of
\cite[Lemma 8.7]{SS} to show that the dual of the cokernel of $\al$
is a finitely generated torsion $\Zp\ps{H/N}$-module, and therefore,
is psuedo-null over $\Zp\ps{G/N}$. Hence we have
 \[ q\big(X(E/F_{\infty})_N\big) = q\big(X(E/L_{\infty})\big), \]
where $q$ is the quotient functor from the category of finitely
generated $\Zp\ps{G/N}$-modules to the category of finitely
generated $\Zp\ps{G/N}$-modules modulo pseudo-null
$\Zp\ps{G/N}$-modules. Since $X(E/L_{\infty})$ is completely
faithful over $\Zp\ps{G/N}$, it follows that $X(E/F_{\infty})_N$ is
a faithful $\Zp\ps{G/N}$-module, as required. \epf

\br \label{control remark on complete faithful}
 It is clear from the proof that under condition (iv)(b), one
 actually shows that $X(E/F_{\infty})_N$ is
a completely faithful $\Zp\ps{G/N}$-module. However, due to the lack
of an analogous result for completely faithful modules in the
direction of Lemma \ref{main lemma}, we are not able to deduce
complete faithfulness of $X(E/F_{\infty})$ from the complete
faithfulness of $X(E/F_{\infty})_N$. This is also precisely the
reason why we require the extra assumption that $r=1$ in condition
(iv)(b). \er

In the next proposition, we mention the best we can do when we do
not assume $r=1$ in condition (iv)(b) which might be of interest.

\bp
 Retaining the assumptions $(i)$, $(ii)$ and $(iii)$
of Proposition \ref{faithful control}. Furthermore, we assume that
the action of $G$ on $N_i/N_{i-1}$ by inner automorphism is given by
a homomorphism $\chi_i:G/N\lra \Zp^{\times}$ for every $i$. Suppose
that for every $v\in S$, the decomposition group of
$\Gal(L_{\infty}/F)$ at $v$ has dimension $\geq 2$, and suppose that
$X(E/L_{\infty})$ is a completely faithful $\Zp\ps{G/N}$-module.
Then $X(E/F_{\infty})_N$ is a completely faithful
$\Zp\ps{G/N}$-module, and for $i\geq 1$, $H_i(N, X(E/F_{\infty}))$
is either a pseudo-null $\Zp\ps{G/N}$-module or a completely
faithful $\Zp\ps{G/N}$-module. \ep

\bpf
 The proof of Proposition \ref{faithful control} carries over to
show that $X(E/F_{\infty})_N$ is a completely faithful
$\Zp\ps{G/N}$-module. By \cite[Proposition 4.2]{Ka06},
$H_i(N,X(E/F_{\infty}))$ is a successive extension of twists of a
$\Zp\ps{G/N}$-subquotient $T$ of $X(E/F_{\infty})_N$ by a one
dimensional character. Therefore, if $H_i(N, X(E/F_{\infty}))$ is
not a pseudo-null $\Zp\ps{G/N}$-module, then $T$ cannot be a
pseudo-null $\Zp\ps{G/N}$-module. Since $X(E/F_{\infty})_N$ is
completely faithful over $\Zp\ps{G/N}$, so is $T$. It is not
difficult to verify that every twist of $T$ by a one dimensional
character is also completely faithful over $\Zp\ps{G/N}$. Hence we
may apply Lemma \ref{faithful ext} to conclude that $H_i(N,
X(E/F_{\infty}))$ is completely faithful over $\Zp\ps{G/N}$. \epf

We also mention that it is clear from the proof of Proposition
\ref{faithful control} that we can prove the following proposition
for a general $N$.

\bp \label{faithful control descent}
 Let $E$ be an elliptic curve over $F$
with either good ordinary reduction or multiplicative reduction at
every prime of $F$ above $p$.
 Let $F_{\infty}$ be a strongly admissible $p$-adic Lie extension of $F$ with
 $G= \Gal(F_{\infty}/F)$.
  Suppose that the following statements hold.
  \begin{enumerate}
\item[$(i)$] $N$ is a closed normal subgroup of $G$ which is contained
in $H$.
\item[$(ii)$] $G/N$ is a non-abelian pro-$p$ $p$-adic Lie group without
$p$-torsion. $($In particular, the dimension of the $p$-adic Lie
group $G/N$ is necessarily $\geq 2$.$)$
\item[$(iii)$] Set $L_{\infty} := F_{\infty}^N$. Suppose that either $(a)$ or $(b)$ holds.
 \begin{enumerate}
 \item[$(a)$] $E_{p^{\infty}}$ is not rational over $F_{\infty}$ $($note the slight difference
 here$)$ and
 $X(E/L_{\infty})$ is a faithful $\Zp\ps{G/N}$-module.
 \item[$(b)$] For every $v\in S$, the decomposition group of $\Gal(L_{\infty}/F)$
  at $v$ has dimension $\geq 2$, and
  $X(E/L_{\infty})$ is a completely faithful $\Zp\ps{G/N}$-module.
  \end{enumerate}
  \end{enumerate}
 Then $X(E/F_{\infty})_N$ is a faithful $\Zp\ps{G/N}$-module.
\ep

The next control theorem is in the direction of a Hida deformation.
We recall that $A$ is the $R$-cofree Galois module attached to the
Hida deformation as defined at the end of Section 4, where
$R=\Zp\ps{X}$, and has the property that $A[P] = E_{p^{\infty}}$ for
some prime ideal $P$ of $R$. As before, we denote by
$X(A/F_{\infty})$ the dual Selmer group of the Hida deformation.

\bp \label{faithful control for Hida}
 Let $F_{\infty}$ be a strongly admissible $p$-adic Lie extension of $F$
 with Galois group
 $G$. Suppose that the following statements hold.
  \begin{enumerate}
\item[$(i)$] $G$ is non-abelian and has dimension $\geq 2$.
\item[$(ii)$] $X(E/F_{\infty})$ is torsion over $\Zp\ps{G}$.
\item[$(iii)$] Either $(a)$ or $(b)$ holds.
 \begin{enumerate}
 \item[$(a)$] $E_{p^{\infty}}$ is not rational over $F_{\infty}$ and
 $X(E/F_{\infty})$ is a faithful $\Zp\ps{G}$-module.
 \item[$(b)$] For every $v\in S$, the decomposition group of $G$
  at $v$ has dimension $\geq 2$, and
  $X(E/F_{\infty})$ is a completely faithful $\Zp\ps{G}$-module.
\end{enumerate}
\end{enumerate}
 Then $X(A/F_{\infty})$
 is faithful over $R\ps{G}$.
\ep

\bpf
 The proof is essentially similar to that in Proposition \ref{faithful
 control}. The only thing which perhaps requires additional attention
 is to show that the cokernel of the map
  \[ \beta: X(A/F_{\infty})/P \lra X(E/F_{\infty}) \]
is finite under the assumption of condition (iii)(a). Note that the
dual of its cokernel is contained in $A(F_{\infty})/P$, where we
write $A(F_{\infty}) = A^{\Gal(\bar{F}/F_{\infty})}$. Noting that
$A[P] = E_{p^{\infty}}$, it then follows from (the dual of) Lemma
\ref{torsion x} that
\[ \corank_{\Zp} E(F_{\infty})_{p^{\infty}} =\corank_R A(F_{\infty})
+\corank_{\Zp}A(F_{\infty})/P.\]
 Again, by the first assumption of condition (iii)(a) and
\cite[Proposition 10]{Z04}, we have that
$E(F_{\infty})_{p^{\infty}}$ is finite. Combining this observation
with the above equation, we have that $\corank_{\Zp} A(F_{\infty})/P
= 0$, or equivalently, that $A(F_{\infty})/P$ is finite. This in
turns implies that the cokernel of $\beta$ is finite. The remainder
of the proof proceeds as in Proposition \ref{faithful control}. \epf

\br
 Of course, one can have a control theorem result for the faithfulness of
 the Selmer group of the Hida
 deformation for other specializations. In particular, one can also
 generalize the above control theorem to (appropriate) Selmer groups
 of more general deformations
 over $\Zp\ps{X_1,..., X_r}$ and their various intermediate
 specializations as considered in \cite{Gr}.
\er

We end the paper discussing an example to illustrate our control
theorem results. Let $p = 5$. Let $E$ be the elliptic curve $21a4$
of Cremona's tables given by
 \[ y^2 + xy = x^3 + x\]
and let $A$ be an elliptic curve $1950y1$ of Cremona's tables
 \[A: y^2 + xy = x^3 - 355303x  - 89334583. \]
Let $p= 5$, $F=\Q(\mu_5)$ and $L_{\infty} = F(A_{5^{\infty}})$. As
discussed in \cite[Section 7]{BZ}, if $X(E/F)$ is finite (as
suggested by its $p$-adic $L$-function), then $X(E/L_{\infty})$ is a
completely faithful $\Z_5\ps{\Gal(L_{\infty}/F)}$-module. \textit{We
will assume this latter property throughout our discussion here}.
Let $F_{\infty}$ be any strongly admissible $5$-adic Lie extension
of $F$ which contains $L_{\infty}$ and such that $N =
\Gal(F_{\infty}/L_{\infty})$ satisfies the condition in Proposition
\ref{faithful control}. It is not difficult to see that
$E_{5^{\infty}}$ is not rational over $F_{\infty}$. Hence we can
apply Proposition \ref{faithful control} to conclude that
$X(E/F_{\infty})$ is a faithful
$\Z_5\ps{\Gal(F_{\infty}/F)}$-module. Examples of strongly
admissible $5$-adic extensions $F_{\infty}$ that we may take are:
\[\Q(A[5^{\infty}], 2^{5^{-\infty}}),\quad
\Q(A[5^{\infty}], 2^{5^{-\infty}}, 3^{5^{-\infty}}),
\quad\Q(A[5^{\infty}], 3^{5^{-\infty}}, 5^{5^{-\infty}}),\quad
\Q(A[5^{\infty}], 2^{5^{-\infty}}, 3^{5^{-\infty}},
5^{5^{-\infty}})\]
\[ \Q(A[5^{\infty}], 2^{5^{-\infty}}, 3^{5^{-\infty}},
5^{5^{-\infty}}, 7^{5^{-\infty}}), \quad \Q(A[5^{\infty}],
2^{5^{-\infty}}, 3^{5^{-\infty}}, 5^{5^{-\infty}}, 7^{5^{-\infty}},
11^{5^{-\infty}}),~...~\mbox{etc},\]
\[ M_{\infty}(A[5^{\infty}],
2^{5^{-\infty}}), \quad M_{\infty}(A[5^{\infty}], 2^{5^{-\infty}},
3^{5^{-\infty}}), \quad M_{\infty}(A[5^{\infty}], 3^{5^{-\infty}},
5^{5^{-\infty}}),\] \[M_{\infty}(A[5^{\infty}], 2^{5^{-\infty}},
3^{5^{-\infty}}, 5^{5^{-\infty}}) \quad M_{\infty}(A[5^{\infty}],
2^{5^{-\infty}}, 3^{5^{-\infty}}, 5^{5^{-\infty}},
7^{5^{-\infty}}),~...~\mbox{etc}, \] where $M_{\infty}$ is any
$\Z_5^r$-extension of $F$ disjoint from $F^{\cyc}$ for $1\leq r\leq
2$. However, at present, we are not able to determine whether or not
$X(E/F_{\infty})$ is completely faithful for any of such
$F_{\infty}$.

\end{document}